\documentclass[11pt]{amsart}

\title{Remarks on a Problem of Eisenstein}
\author{Roger C. Alperin}
\address{
Department of Mathematics\\  San Jose State University\\ San Jose,
CA 95192, USA}
\email{alperin@math.sjsu.edu}
\usepackage{amsmath}
\usepackage{amsthm}
\usepackage{amscd}
\usepackage{amssymb}
\usepackage{verbatim}

\newcommand{\A}{\mathbb A}

\newcommand{\Q}{\mathbb Q}
\newcommand{\Z}{\mathbb Z}

\newcommand{\e}{\epsilon}

\newtheorem{theorem}{Theorem}[section]
\newtheorem{proposition}[theorem]{Proposition}
\newtheorem{lemma}[theorem]{Lemma}

\newtheorem{example}[theorem]{Example}

\def\o{\mathcal O}
\def\A{\mathcal A}

\def\n{\sqrt{N}}

\def\endproof{~$\blacksquare$\medskip}

\date{}
\begin{document}

\maketitle

\begin{abstract} The fundamental
unit of $\Z[\sqrt{N}]$ for square-free $N=5\ mod\ 8$ is either $\epsilon$
or $\epsilon^3$ where $\epsilon$ denotes the fundamental unit of the maximal
order of $\Q(\sqrt{N})$. We give infinitely many examples for each case.
\end{abstract}

\section{Introduction} 
For $N$ square-free, the ring of integers $\o_N$ of a
real quadratic field
$\Q(\sqrt{N})$ 
has an infinite cyclic group of units of index 2.
The  generator $\e$ for this subgroup is 
the fundamental unit.
The ring of integers
$\o_N$ has a subring $\A_N=\Z[\sqrt{N}]$; this is a proper subring if and only
if
$N=1\mod 4$. The subring also has an infinite cyclic subgroup of units
generated by
$\e^e$; it is easy to see that $e=1$ or $e=3$; the latter  occurs only if 
$N=5 \mod 8$. 

Characterizing those $N$ for which $e=3$ is 
the problem of Eisenstein in the title of this article. 
By elementary methods we shall give
infinitely many examples  for each of the cases of
$e=1$ or
$e=3$. This problem has 
been addressed in \cite{IW} and \cite{S} using other methods.

\section{Main Examples}
Basic properties of continued fractions and the relation of equivalence can be
found in \cite{HW}.  Equivalence of two continued fractions means that the
periodic parts  are equal or equivalently that the two real numbers are
related by a linear fractional transformation.

The following examples are well-known
\cite[p. 297]{S}:

\begin{example}
$\sqrt{a^2+4}=(a;\overline{\frac{a-1}{2},1,1,\frac{a-1}{2},2a})$
%=
%a+\cfrac{1}{{\frac{a-1}{2}}+\cfrac{1}{{1}+\cfrac{1}{{1}+
%\cfrac{1}{{\frac{a-1}{2}}+\cfrac{1}{{2a}+\dotsb}}}}}$
for any odd integer
$a>1$.
\end{example}

Consider $a=4b\mp 1$ and $N=a^2+4$ then 
\begin{equation*}
\frac{1}{\frac{\n\pm 1}{4}- b}=\frac{4}{\n-a}\frac{\n+a}{\n+a}
=4\frac{\n+a}{N-a^2}=\n+a
\end{equation*}

\begin{proposition} Suppose $a$ is odd and greater than 1. For $N=a^2+4$, 
then 
$\frac{\n \pm 1}{4}$ is equivalent to   $\n$.
\end{proposition}
\proof For $a=4b\mp 1$ the floor of $\frac{\n \pm 1}{4}$ is $b$. 
\endproof

\begin{example}
For any odd integer
$a>3$,
$\sqrt{a^2-4}=(a-1;\overline{1,\frac{a-3}{2},2,\frac{a-3}{2},1,2a-2})$. 
%$\sqrt{a^2-4}=a-1+\cfrac{1}{{1}+\cfrac{1}{{\frac{a-3}{2}}+\cfrac{1}{{2}+
%\cfrac{1}{{\frac{a-3}{2}}+\cfrac{1}{{2a-2}+\dotsb}}}}}$.
\end{example}

As a consequence one can easily show that 
$$1+\frac{\sqrt{a^2-4}}{a-2}=(2;\overline{\frac{a-3}{2},1,2a-2,1,\frac{a-3}{2}}).$$

Let $N=a^2-4$ and put $a=4b\pm1$.
For $a=4b-1$ we have
\begin{equation*}
\frac{1}{\frac{\n - 1}{4}-
(b-1)}=\frac{4}{ {\n} -(a-2)}=\frac{ {\n} +(a-2)}{a-2}.
\end{equation*}
For $a=4b+1$ we obtain
\begin{equation*}
\frac{1}{\frac{\n + 1}{4}-
b}=\frac{4}{ {\n} -(a-2)}=\frac{ {\n} +(a-2)}{a-2}.
\end{equation*}

\begin{proposition} Suppose $a$ is odd and greater than 3. For $N=a^2-4$
 then 
$\frac{\n \pm 1}{4}$ is equivalent to   $\n$.
\end{proposition}
\proof For $a=4b\pm 1$ we have $\frac{\n \pm 1}{4}$ is equivalent to
$1+\frac{\n}{a-2}$ which is equivalent to $\n$. 
\endproof

\begin{example}
For any integer
$a>1$ $\sqrt{a^2+1}=(a;\overline{2a})$.
%=a+\cfrac{1}{{2a}+\cfrac{1}{{2a}+\dotsb}}$ .
\end{example}

\begin{proposition} \label{length}For $N=4a^2+1$ where $a$ is odd and greater than
3, then 
$\frac{\n \pm 1}{4}$ is not equivalent to   $\n$.
\end{proposition}
\proof The numbers $u_{\pm}=({\frac{\n \pm 1}{4}-\lfloor{\frac{\n \pm
1}{4}}\rfloor})^{-1}$ are greater than 1 by
definition.  They are purely periodic (\cite{HW}) since the conjugates are 
negative and 
$-\frac{1}{\bar{u}_{\pm}}={\frac{\n \mp 1}{4}+\lfloor{\frac{\n
\pm 1}{4}}\rfloor}$ 
is greater than 1.  

If
$\frac{\n
\pm 1}{4}$ is equivalent to
$\n$ then $u_{\pm}$ has period length one also. Hence $u_{\pm}=(\overline{2a};)$.
The continued fraction $(\overline{2a};)$ satisfies the equation $x^2-2ax-1$ 
which has the solutions 
$\sqrt{a^2+1}\pm a$; these can not be the same as $u_{\pm}$. This
contradiction gives the desired result.
\endproof

\section{Relations of Units to Continued Fractions}
We suppose that $N=5 \mod 8$ is square-free.  It is an elementary exercise to
see that the fundamental unit
$\e$ is a solution to $x^2-Ny^2=\pm 4$ with $x$, $y$ odd if and only if 
$e=3$. 

Let
$\A=\A_N$ and $\o=\o_N$. Consider the ideals $I_{\pm}=[4,\n\pm 1]$ in $\A$.
(the generators are a lattice basis).
Extend these ideals to ideals $J_{\pm}=2[2,\frac{\n\pm 1}{2}]$ in $\o$; thus
$J_{\pm}$ is  principal  since when $N=5 \mod 8$ the ideal
(2) is   maximal. An easy calculation shows that $[4,\n+1]^2=2[4,\n-1]$ so
that 
 $[4,\n+1]$ is an element of order 1 or 3 in the class group $Cl(\A)$.

\begin{lemma} \label{key} When $N=5 \mod 8$ the following are equivalent:

(a) The equation $x^2-Ny^2=\pm 4$ has a solution with odd integers $x, y$.

(b) There is a non-integral element of norm $\pm 4$ in $\A_N$.

(c) The ideals $I_{\pm}$ are principal.

(d) The elements $\frac{\n\pm 1}{4}$ are equivalent to $\n$.
\end{lemma}
\proof It is easy to see that (a) and (b) are equivalent using $N=5 \mod 8$.
The conditions (b) and (c) are also easily seen to be equivalent since the
ideals
$I_{\pm}$ have norm 4. Conditions (c) and (d) are equivalent using the
well-known description of the class group in terms of equivalence classes of
elements according to their continued fractions.
\endproof

If the elements $\frac{\n\pm 1}{4}$ are not on
the principal cycle then the two continued fractions are the reverse of one
another since the elements $[4, \n \pm 1]$ are inverses of one another in
the class group of $\A$. 

\begin{theorem} Suppose $N=5 \mod 8$ is square-free. Consider the
surjective natural homomorphism $$\phi: Cl(\A_N)\rightarrow Cl(\o_N).$$

(a) The homomomorphism $\phi$ is an isomorphism if and only if
$e=3$.

(b) The homomorphism $\phi$ has kernel
generated by $[4,\n+1]$ if and only if $e=1$.
\end{theorem}
\proof It is well-known that $\phi$ is surjective,  that the kernel
has order dividing three, and the order of the kernel is three if and only if
condition (a) of the Lemma fails (\cite{St}). 
 Using  Lemma
\ref{key} and this remark we see that the kernel of $\phi$ is 
the ideal class of $[4,
\n+1]$, and hence this class is an element of of order 3 if and only if $e=1$. 
\endproof

\section{Applications}
Using a theorem of Erd\u{o}s
(\cite{E}) it follows that there
are infinitely many square-free integers  $a^2\pm 4$ or $4a^2+1$ for odd $a$ . 

\begin{theorem} For $a$ odd and greater than 3. There are infinitely many 
square-free
$N=4a^2+1$ with $e=1$.
\end{theorem}
\proof It follows from Proposition \ref{length} that $\frac{\n\pm
1}{4}$ have cycle lengths greater 1 and hence are not equivalent to $\n$; thus
the ideals 
$[4, \n\mp 1]$ of
$\A_N$ are not
 principal and therefore there is no element of norm 4 so the fundamental unit
$\e$ does  belong to
$\A_N$; hence
$e=1$.
\endproof

\begin{theorem} For $a$ odd and greater than 3.  There are infinitely many
square-free $N=a^2\pm 4$ with $e=3$.
\end{theorem}
\proof
The numbers
$u_{\pm}=\frac{\n
\pm 1}{4}$ are  equivalent to $\n$.  Consequently the ideal $[4, \n\mp 1]$ of
$\A_N$ is
 principal and therefore the fundamental unit $\e$ does not belong to $\A_N$;
hence
$e=3$.
\endproof

\smallskip
\par\small{ MSC: 11R65,11R29}

\begin{thebibliography}{ABC}
\bibitem{E} P. Erd\u{o}s, {\it Arithmetical propoertes of polynomials}, J. London
Math. Soc., 28, 1953, 416-425.
\bibitem{HW} G. H. Hardy and E. M. Wright, {\bf An Introduction to the Theory of
Numbers}, Oxford, Fifth Edition, 1979.
\bibitem{IW} N. Ishii, P. Kaplan, K. S. Williams, {\it On Eisenstein's Problem},
Acta Arithmetica,  54, 1990, 323-345.
\bibitem{S} W. Sierpinski, {\bf Elementary Theory of Numbers}, Polish Academy of
Sciences, Vol. 42, Warsaw, 1964.
\bibitem{St} P. Stevenhagen, {\it On a Problem of Eisenstein}, Acta Arithmetica,
74, 3, 1995, 259-26.
\end{thebibliography}
\end{document}